\documentclass[journal]{IEEEtran}
\usepackage{graphicx}
\usepackage{algorithm2e}
\usepackage{amsmath}
\usepackage{amssymb}
\usepackage{url}
\usepackage{amsmath,amsfonts,amsthm}
\usepackage{array}
\usepackage[caption=false,font=normalsize,labelfont=sf,textfont=sf]{subfig}
\usepackage{textcomp}
\usepackage{stfloats}
\usepackage{url}
\usepackage{verbatim}
\usepackage{cite}

\newtheorem{thm}{Theorem}[section]
\newtheorem{lem}[thm]{Lemma}

\theoremstyle{definition}

\DeclareMathOperator*{\argmin}{arg\,min}

\begin{document}

\title{A direct solution to the interpolative inverse non-uniform fast Fourier transform (i-iNFFT) problem for spectral analyses of non-equidistant time-series data}

\author{Michael~Sorochan~Armstrong,
        Jos\'{e} Carlos P\'{e}rez-Gir\'{o}n,
        Jos\'{e} Camacho \IEEEmembership{Member,~IEEE},
        Regino Zamora

\thanks{\IEEEcompsocthanksitem Michael Sorochan Armstrong and Jos\'{e} Camacho are part of the Computational Data Science (CoDaS) Lab in the Department of Signal Theory, Telematics, and Communications at the University of Granada.\protect\\

E-mail: mdarmstr (at) go.ugr.es

\IEEEcompsocthanksitem Jos\'{e} Carlos P\'{e}rez-Gir\'{o}n and Regino Zamora are part of the Inter-university Institute for Research on the Earth System in Andalucia through the University of Granada.

}
}
\markboth{IEEE Transactions on Signal Processing}%
{Shell \MakeLowercase{\textit{et al.}}: Bare Demo of IEEEtran.cls for Computer Society Journals}

\maketitle

\IEEEtitleabstractindextext{%
\begin{abstract}
Recovering the most significant oscillatory patterns in times-series climate measurements using Fourier analyses assumes a constant sampling rate, but long periods of sensor inactivity violate this assumption which lead to inaccurate results. This is obvious when the calculated spectrum of Fourier coefficients is not sufficient to reconstruct the original values of the time series via the inverse Fourier transform.
This work demonstrates a closed-form solution to the Interpolative Inverse Non-Uniform Fast Fourier Transform (i-iNFFT) which can be used to reconstruct the original values in a time series and interpolate missing values, in theoretical $\mathcal{O}(MlogM + N^2)$ time complexity for $M$ observations in the time domain and $N$ Fourier coefficients in the frequency domain. The proposed solution's theoretical time complexity is competitive with the most current research on the topic of Non-Uniform Fast Fourier Transforms (NFFTs) and their inversions \cite{kircheis2019direct}. Our results suggest that for missing, contiguous blocks of time-series data the i-iNFFT better predicts procedurally omitted data when compared with an equivalent truncated, inverse Fast Fourier Transform (t-iFFT) for an equal number of Fourier coefficients that assumes a constant sampling rate. These results can be further improved using stochastic gradient descent to infer the relative scale of a periodic function from an unlabelled time-series dataset.
\end{abstract}

}

\IEEEdisplaynontitleabstractindextext

\section{Introduction}
\IEEEPARstart{M}{uch} of experimental science has experienced a paradigm shift towards collecting larger volumes of data to develop and test more complicated scientific hypotheses \cite{sanchez2018integrated,grabherr2005long}. For environmental monitoring, these data are often collected in inhospitable environments where human curation of the collection devices is limited. This lends to time series measurements with significant periods of missing or unreliable measurements due to disruptions in network connectivity, power supply, or adverse environmental conditions. But the use of remote sensor data is indispensable for making informed decisions on natural resource management and environmental protection - particularly in hard-to-reach ecosystems such as high mountain environments. This type of meteorological information frequently lacks adequate coverage due to its origin far from stations located at lower altitudes and nearer to urban areas where they can be more easily accessed.


Binning time-series data can decrease the impact of missing values on a resultant model by measuring the average of several ``binned'' observations via a functional reduction in sampling rate. However binning does not make use of the full acquisition rate of the instrument, so this process negates the benefits of high-throughput measurements. Even by way of a moving average, in cases where the number of observations within the bin decrease towards the size of the expected measurement gaps it becomes impractical or impossible to reasonably infer the missing values.

Binning is neither ideal, nor universally applicable in the case of irregularly sampled time-series data. A broader representation of the data can be formulated via an analysis of its corresponding \textit{spectrum} of frequency-domain components, whose complex-valued coefficients for an equi-dimensional transformation encode the relative offset and amplitude of a series of additive sinusoidal functions that best reconstruct the original time-series data. But the Discrete Fast Fourier Transform (DFT), and its corresponding inverse that are most commonly used to perform this type of analysis assume a sampling rate that is constant on average \cite{friel1999frequency}. 

Least-squares spectral analysis of irregularly sampled data is an active area of research that commonly finds applications wherever there is an interruption in measurement - as is commonly seen in climate research \cite{haley2021missing}, and other environmental monitoring \cite{dilmaghani2007harmonic}. To minimize the effect of spectral leakage, an apodization (i.e. weight) function can be applied to discount the contributions of higher-order variability as a low-pass filter. The characteristics of this filter can be derived from \textit{a-priori} mathematical knowledge, to reduce user intervention and influence in the analysis \cite{potts2021continuous}.


Interpolative Inverse Non-uniform Fast Fourier Transforms (i-iNFFTs) calculate a series of Fourier coefficients such that the original data, presenting with irregular sampling, can be reconstructed in a way that interpolates a reasonable amount of measurement gaps via the preferential weighting of lower-frequency components. As a series of complex numbers, Fourier coefficients are simple descriptors with real physical properties that can be interpreted and analyzed for how they relate to experimental data. But previous implementations\cite{kunis2007stability} of this routine do not exploit a closed-form solution and convergence is difficult to measure - owing to a reliance on minimizing numerical error rather than residuals (i.e. a similarity to the ``true'' solution). The original implementation of the i-iNFFT algorithm must also be modified at a low level according to the specific needs of the analysis before being called from a higher-level 
routines \cite{KeKuPo09}, making it more difficult to use for many applications. This work presents:
\begin{enumerate}
    \item a mathematical description of the limitations of the Discrete Fast Fourier Transform (DFT) for irregularly sampled data (\textbf{Appendices \ref{sec:app1} \ref{sec:app2}}), and an illustration of the possible variants of NFFTs and their inversions (\textbf{Section \ref{sec:dft}}).
    \item Previous work on the topic (\textbf{Section \ref{sec:lssa}}).
    \item Derivation and proof of a closed-form solution to the optimal interpolation problem (\textbf{Section \ref{sec:solve}}), and an overview of the criteria for selecting a weight function (\textbf{Section \ref{subsec:weights}}).
    \item An analysis of the internal and external reconstruction errors as a function of the number of Fourier coefficients, and the relative weights ascribed to them (\textbf{Section \ref{sec:res}})
    \item A method for inferring optimal labels for time, when the precise periodicity is unknown (\textbf{Section \ref{subsec:time_labels}}).
    \item A highly portable Python package using standard libraries.
\end{enumerate}
\subsection{Notation conventions}
We define a continuous series as $f_j$ and $h_k$ for $j \in [0,M-1]$ in the time domain and $k \in [-N/2, N/2-1]$ in the frequency domain. The \textit{hat} designator indicates a series of predicted values in cases where the transformation is appropriately scaled such that its inverse is known, and in cases where it is not. Matrices are denoted as capital letters, except $M$ and $N$ which refer to the number of samples in the time domain and the number of Fourier coefficients. In instances where either the time or frequency series are discontinuous, notation of the form $f(x_j), h(x_k), ...$ etc. are used to indicate this. Frequency-domain weights are denoted as $w_k$, and $W_k$ for the vector and diagonal matrix forms. The weighted norm, denoted using $||x||_{W^{-1}_k}$ for some complex vector $x$ is evaluated as $x^HW^{-1}_kx$.

\subsection{Intuitions for the inverse transform}
The forward and inverse Fourier transforms on regularly-spaced data can be thought of as linear transformations from the time to frequency domains, and frequency to time domains respectively. The dimensions of this transformation matrix are related to the number of observations in the time domain $M$, and the $N$ number of expected Fourier coefficients centred about the limiting Nyquist frequency, $\pm N/2$. For the equidistant case, the Hermitian (i.e. the conjugate transpose) of the transformation matrix $F \in N \times M$ is proportional to its inverse, that is: $F^H \propto F^{-1}$ via a scaling factor $1/M$. $F$ is an \textit{orthogonal} (see: Appendix \ref{sec:app1}), and \textit{Vandermonde} matrix that allows for a efficient and reversible transformation via the (discrete) Fast Fourier Transformation (FFT) algorithm \cite{petersen2008matrix}.

If the data is irregular (i.e. non-equispaced) in either the frequency or time domain, the dimensionality of the transformation matrix $A \in N \times M$ is still defined by the number of time domain observations and anticipated frequency components in the case of the forward transform. However due to the irregular characteristics of either the time or frequency components encoded by the transformation matrix, it is no longer sufficient to simply take its Hermitian as the inverse, since the matrix is neither orthogonal (in such a way that the inverse transform can be appropriately scaled by a scalar factor, see: Appendix \ref{sec:app1}), nor Vandermonde (such that an efficient algorithm can be used to perform the transformation). That is to say: $A^H \not\propto A^{-1}$.

For large values of $M$ and $N$, it is not practical to directly invert the forward transformation matrix for the inverse transform. Rather, we discuss instead whether the vector in the uni-dimensional case \textit{predicted} by the Forward transform is \textit{invertible} through the application of the transformation matrix's Hermitian on the predicted values themselves. In this way, if the predicted vector of Fourier coefficients faithfully recover the original signal, then we say that the coefficients themselves are invertible, rather than discussing the inversion of the transformation matrix itself as demonstrated in Section \ref{subsec:inv}. Previous work has focused on the fact that it is more efficient to calculate invertible coefficients than it is to directly invert the transformation matrix as is easily done in the equidistant case. Therefore the calculation of the inverse performs a similar operation from the time to frequency domain, except the Fourier coefficients are calculated in such a way that the original signal is recoverable without any additional scaling information as is available in the equidistant case demonstrated in Appendix \ref{sec:app1}.

The rank of the irregular transformation matrix $A$ was shown in previous work to be related to the minimum distance spanned by two measurements, $q$ and the $d^{th}$-order (i.e. dimensionality) of the transform. $A$ has been shown to \textit{under-determined}, when the number of frequency coefficients $N > 2dq^{-1}$, and \textit{over-determined} when $N < 2dq^{-1}$ \cite{kunis2007stability}. 



\subsection{Interpolative inverse transforms}
An interpolative transform is any transformation from the frequency domain to the time domain such that the signals not populated by time-domain measurements are estimated from the Fourier coefficients. This can be done in order to up-sample otherwise regular time-domain measurements \cite{oppenheim1999discrete}, but an additional use case arises where there are relatively large sections of missing data from an otherwise regularly measured time-series. Interpolation is done by weighting the lower-frequency contributions of the Frequency spectrum more strongly according to some criteria.

The Inverse Non-Uniform Fast Fourier Transform (iNFFT) calculates a spectrum of frequency coefficients that are invertible through the Hermitian of the forward transformation matrix. An \textit{interpolative} inverse Non-Uniform Fast Fourier Transform (i-iNFFT) includes an additional constraint which restricts the solution via a selected apodization function so that subsequent transformation back from the frequency to time domain optimally interpolates the missing data. Use of an apodization function improves the numerical characteristics of the transformation, allowing the user to select a higher number of Fourier coefficients for the reconstruction of the time-series data than they would typically need. An overview of the criteria for selecting an apodization function is presented in Section \ref{subsec:weights}, following the work of Kunis and Potts \cite{kunis2007stability}.

\section{Fourier analysis}\label{sec:dft}
\subsection{Discrete Fourier Analysis}
The Discrete Fourier Transform, (DFT) and the related Fast Transforms (FFTs) \cite{cooley1965algorithm} assume that the time domain data are of consistent density ($j/M$). In the case of a presumably continuous time series signal, $f_j \in \mathbf{R}$ where $j/M \in [-\frac{1}{2},\frac{1}{2})$, the DFT is calculated by evaluating the time-domain series to get the complex Fourier coefficients $h_k \in [-\frac{N}{2}, \frac{N}{2})$, where $h_k \in \mathbf{C}$:
\begin{equation}\label{eq:dft}
    h_k = \sum_{j=-M/2}^{M/2-1} f_j e^{-2\pi \mathbf{i} k \frac{j}{M}}
\end{equation}
\noindent and its inversion evaluates $h_k$ to yield the Fourier coefficients $f_j$:
\begin{equation}\label{eq:idft}
    f_j = \frac{1}{M}\sum_{k=-N/2}^{N/2 - 1} h_k e^{2\pi \mathbf{i} \frac{j}{M} k}.
\end{equation}
When the number of Fourier coefficients is fewer than the number of observations, (i.e. $N < M$), the equality in Eq. (\ref{eq:idft}) holds only up to an approximation, and is referred to as a \textit{truncated} Discrete Fourier transfrom. The full spectrum of Fourier coefficients in $N$ can be calculated up to $M = N$ where the maximum frequency of $\pm N/2$ is the Nyquist frequency, beyond which the sampling rate is not sufficient to warrant calculation of any further coefficients. A full Fourier representation is theoretically sufficient to fully reconstruct the original data, $f_j$ - however there are some notable exceptions in cases of discontinuous data that are not readily represented by a series of trigonometric polynomials \cite{carslaw1925historical}. It is also important to note that the convention is sometimes reversed \cite{rasheed2009fourier,van1992computational} between the \textit{forward} and \textit{inverse} transforms by the signage of the exponential term and normalization coefficient in Eqs. (\ref{eq:dft}) and (\ref{eq:idft}). 

These transformations can be written as a linear operation $f_k = Ff_j$, where Eqs. (\ref{eq:dft}) and (\ref{eq:idft}) become the transformation matrices:
\begin{align}\label{eq:linear_dft}
    F = \left(e^{-2\pi\mathbf{i}k\frac{j}{M}}\right)_{j,k \in N,M}\\
    F^{-1} = \frac{1}{M}\left(e^{2\pi\mathbf{i}\frac{j}{M}k}\right)_{j,k \in N,M}.
\end{align}
Since these matrices are Vandermonde, both the forward operation: $h_k = Ff_j$ and the inverse: $f_j = F^{-1}h_k$ can be performed in $\mathcal{O} (M log M)$ time complexity using the FFT algorithm, which recursively divides the series into smaller operations (assuming constant spectral density via the $j/M$ term). Even in cases where the input data are continuous, but a fewer number of Fourier coefficients are calculated, the truncated Fourier transform performs the operation in $\mathcal{O} (N log M)$ time \cite{markel1971fft} \cite{bowman2018partial}.
A fundamental property of the DFT/FFT is its self-inversion property:
\begin{equation}\label{eq:self1}
    \frac{1}{M}F^H = F^{-1},
\end{equation}
and so:
\begin{equation}\label{eq:self2}
FF^H = MI_N. 
\end{equation}
\noindent This identity implies orthogonality of the transformation matrix $F$.
\subsection{Non-Uniform Fast Fourier Transforms}
It is possible to adapt the linear transformations in Eq. (\ref{eq:linear_dft}) for irregularly sampled data in the time domain: $f(x_j), x_j/M \in [-\frac{1}{2},\frac{1}{2})$ and/or irregularly sampled data in the frequency domain: $h(x_k), x_k \in [-\frac{N}{2},\frac{N}{2})$ for NDFTs. In terms of their impact on the resultant transformation matrices (here denoted as $A$ for non-equidistant cases), the value of the entries do not increase at a predictable rate as would be the case for the equivalent uniform DFT. The resulting matrices are not Vandermonde, and so the transformation cannot utilize the FFT algorithm directly to accomplish the task in $\mathcal{O}(N log N)$ time \cite{roth1990application}. Exactly how the structure of these matrices differ depends on the class of NDFT they represent - which are related to which domains the measurement irregularities are observed. The most commonly applicable type of NDFT are the \textit{Type-I} transformations Eq. (\ref{eq:typeI}), with irregular sampling in the time, but not frequency domain. \textit{Type-II} Eq. (\ref{eq:typeII}), and \textit{Type-III} Eq. (\ref{eq:typeIII}) transforms assume discontinuities in the frequency but not time, and both time and frequency domains, respectively:
%
%
\begin{equation}\label{eq:typeI}
    A_I = \left(e^{-2\pi\mathbf{i}k\frac{x_j}{M}}\right)_{k,j \in N,M}
\end{equation}
\begin{equation}\label{eq:typeII}
    A_{II} = \left(e^{-2\pi\mathbf{i}x_k\frac{j}{M}}\right)_{k,j \in N,M}
\end{equation}
\begin{equation}\label{eq:typeIII}
    A_{III} = \left(e^{-2\pi\mathbf{i}x_k\frac{x_j}{M}}\right)_{k,j \in N,M}.
\end{equation}
%
%
Non Uniform Fast Fourier Transforms (NFFTs) are algorithmic improvements to NDFTs, analogous to FFTs. NFFTs can perform a Forward transformation in $\mathcal{O}(NlogN + Mlog(1/\epsilon))$ time \cite{keiner2007nfft}, where $\epsilon$ is the desired precision of the result. However in the case of the Non-Uniform Discrete Fourier Transform (NDFT), $AA^H \neq MI_N$, since $A$ is not directly invertible via a constant scaling factor of its Hermitian (Appendix \ref{sec:app2}), and so calculating the inverse requires further consideration.

\subsection{Inverse NFFTs}\label{subsec:inv}
\begin{lem}\label{lm:yes}
In the equidistant case, the inversion of the predicted Fourier coefficients, $\hat{h}_k$, scaled by the inverse of the number of time-domain observations, $1/M$, is invertible by the Hermitian of its forward transformation matrix, $F \in \mathbf{C}^{N\times M}$ because the original continuous time-domain signal $f_j$ is optimally reconstructed up to the precision allowed by the number of Fourier coefficients, $N$. That is:
\begin{equation}
f_j = \frac{1}{M}F^H\hat{h}_k = F^{-1}\hat{h}_k
\end{equation}
\end{lem}
\begin{lem}\label{lm:no}
In the non-equidistant case, the inversion of the predicted Fourier coefficients, $\hat{h}_k$ scaled by the number of time-domain observations, $1/M$, is not invertible by the Hermitian of its forward transformation matrix, $A \in \mathbf{C}^{N \times M}$ and does not reconstruct the original time-domain signal, $f(x_j)$ optimally up to the precision allowed by the number of Fourier coefficients, $N$. That is:
\begin{equation}
    f(x_j) \neq \frac{1}{M}A^H\hat{h}_k \neq A^{-1} \hat{h}_k
\end{equation}
\end{lem}
\begin{thm}\label{thm1}
We define a series of predicted Fourier coefficients $\hat{h}_k, \hat{h}(x_k) : = \hat{h}$ as \textbf{invertible} iff the original signal in the time domain ($f_j, f(x_j) := f$) can be recovered exactly the via Hermitian of its forward transformation matrix (in the general case: $A \in \mathbf{C}^{N \times M}$):
\begin{equation}\label{eq:def1}
f = A^H\hat{h}.
\end{equation}
\end{thm}

Lemmas \ref{lm:yes} and \ref{lm:no} are proven in Appendices \ref{sec:app1} and \ref{sec:app2}. As a consequence of Lemma \ref{lm:no}, in order for the predicted coefficients to satisfy Theorem \ref{thm1}, they must be calculated numerically through the minimization of the following cost function in Theorem \ref{thm2}:

\begin{thm}\label{thm2}
An invertible $\hat{h}_k$ for an irregularly sampled series $f(x_j)$ is one that minimizes the following cost function:
\begin{equation}\label{eq:min_func}
    \argmin_{\hat{h}_k} ||f(x_j) - A^H\hat{h}_k||_F^2.
\end{equation}
\end{thm}
Theorem \ref{thm1} states that in the case of the forward transform $\alpha A^H = A^{-1}$ would need to be satisfied for some constant $\alpha$ if the transform is invertible. For non-equispaced data however $A$ does not have the self-inversion property as demonstrated in Appendix \ref{sec:app2} described in Lemma \ref{lm:no}, so no such $\alpha$ exists. It is reasonable then to modify $\hat{h}_k$ directly via Theorem \ref{thm2} to ``anticipate'' the inversion via an equivalent \textit{adjoint} (equivalent to the \textit{Hermitian} of the transformation matrix previously described) transformation: $A^H$, rather than inverting $A$ directly which would scale in $\mathcal{O}(N^3)$ time. For this reason, we distinguish between ``normalized'' frequency spectra $\hat{h}$ and ``un-normalised'' frequency spectra $\eta$ that do, and do not satisfy Theorem \ref{thm1} and use this to define whether the coefficients themselves are invertible. 
\subsection{Categorizing NFFTs and their inverses}
In the more general case, the NFFT can be thought of relative to either the time or frequency domain depending on which is known \textit{a-priori} to be the ``true'' representation of the data. So NFFTs are categorized based on what domain the irregular sampling takes place in. \textit{Type-I} and \textit{Type-II} are irregularly sampled in the time or frequency domains, while transformations of \textit{Type-III} are irregular in both domains. For NFFTs of \textit{Type-I}, there exist both the un-normalised transforms to and from the time domain, and their respective \textit{inverses} that perform a similar operation, but are distinguished by whether an adjoint transform is anticipated such that the original data is optimally reconstructed as shown in Table \ref{tab:nffts}. 

An adjoint transform from the frequency domain to the time domain can yield either an directly invertible vector of time series data $f$, or $\phi$ depending on whether they are normalized in a way that anticipates a reversed transformation.
\begin{table}[htb!] 
\begin{center}
\begin{tabular}{cc}
\textbf{Standard} & \textbf{Inverse}\\  
\multicolumn{1}{l}{\begin{tabular}[c]{@{}l@{}}Forward NFFT\\ $\hat{\eta}_k = Af(x_j)$\end{tabular}} & \multicolumn{1}{l}{\begin{tabular}[c]{@{}l@{}} Inverse NFFT \\ 
$A^{-1}\hat{h}_k=f(x_j)$\end{tabular}}\\ 
\multicolumn{1}{l}{\begin{tabular}[c]{@{}l@{}}Adjoint NFFT\\ $\hat{\phi}(x_j) = A^Hh_k$\end{tabular}} & \multicolumn{1}{l}{\begin{tabular}[c]{@{}l@{}} Inverse adjoint NFFT \\ $\left(A^{H}\right)^{-1}\hat{f}(x_j) = {h}_k$\end{tabular}} \\
\end{tabular}
\end{center}
\caption{\textbf{Type I NFFTs}: Forward, Inverse, Adjoint, and Inverse Adjoint transforms. Note that for Type-I the time domain is discontinuous and the frequency domain is continuous. The Forward NFFT transforms irregularly-spaced data in the time domain to the frequency domain, but cannot recover the original signal from the adjoint transform. The inverse NFFT performs a similar operation in a numerical fashion such that the original time signal can be recovered from the adjoint transform. An adjoint transform is simply the Hermitian of the forward transformation matrix, and yields un-normalised time-domain data that do not recover the original frequency values via the forward transform. In much the same way, the inverse adjoint transform calculates properly scaled time domain data that can recover the original frequency domain information through the forward transform. Note that $A^{-1}$ and $(A^H)^{-1}$ are not indicative of how the coefficients are calculated, but are indicated using the inverse for illustrative purposes.}
\label{tab:nffts}
\end{table}
For \textit{Type-II} transformations, where irregular sampling is anticipated in the frequency domain - the same considerations hold as shown in Table \ref{tab:nffts2}. A regularly-sampled time-domain function, $f_j$ can be transformed into an irregularly sampled $\hat{\eta}(x_k)$ in cases where an inverse transformation is not anticipated, and therefore $\hat{\eta}(x_k)$ is not invertible as $\hat{h}(x_k)$ would be.
\begin{table}[ht] 
\begin{center}
\begin{tabular}{cc}
\textbf{Standard} & \textbf{Inverse}\\  
\multicolumn{1}{l}{\begin{tabular}[c]{@{}l@{}}Forward NFFT\\ $ \hat{\eta}(x_k) = Af_j$\end{tabular}} & \multicolumn{1}{l}{\begin{tabular}[c]{@{}l@{}} Inverse NFFT \\ 
$A^{-1}\hat{h}(x_k) = f_j$\end{tabular}}\\ 
\multicolumn{1}{l}{\begin{tabular}[c]{@{}l@{}}Adjoint NFFT\\ $ \hat{\phi}_j = A^Hh(x_k)$\end{tabular}} & \multicolumn{1}{l}{\begin{tabular}[c]{@{}l@{}} Inverse adjoint NFFT \\ $\left(A^{H}\right)^{-1}\hat{f}_j = h(x_k)$ \end{tabular}} \\
\end{tabular}
\end{center}
\caption{\textbf{Type II NFFTs}: Forward, Inverse, Adjoint, and Inverse Adjoint. Note thatfor Type-II transforms the frequency domain is discontinuous and the time domain is continuous. Here the Forward transform calculates a series of discontinuous, un-normalized Fourier coefficients that cannot recover the original time domain signal simply via adjoint transform. The inverse transform calculates a similar series of Fourier coefficients such that an adjoint transform would recover the original time series signal. A Type-II adjoint transform calculates the un-normalized, but equidistant time series that cannot recover the original discontinuous Fourier coefficients via the Forward transform. The Inverse Adjoint Transform in this case calculates the continuous time series such that the original discontinuous Fourier coefficients can be recovered using the Forward transform. Note that $A^{-1}$ and $(A^H)^{-1}$ are not indicative of how the coefficients are calculated, but are indicated using the inverse for illustrative purposes.}
\label{tab:nffts2}
\end{table}
%
%
For the purposes of this work and to explicate the far more common scenario where irregularly sampled time-series data are provided, the focus will be exclusively on \textit{Type-I} transforms, although it would be relatively straightforward to apply the same thinking to \textit{Type-II} transforms.
%
\subsection{Optimal interpolation for Fourier-type problems}\label{subsec:weights}

Seeking a continuous function that interpolates the missing values of irregularly sampled dataset, the smallest possible magnitudes of each entry in $\hat{h}_k$ as $N$ Fourier coefficients are calculated as a function of some user-selected weights $w_k$ that also satisfy Eq. \ref{eq:min_func} as a constraint:
\begin{equation}\label{eq:intrp}
    \argmin_{\hat{h}_k | A^H\hat{h}_k = f(x_j)}{\sum_{k \in I_N} \frac{|\hat{h}_k|^2}{w_k}}
\end{equation}
The optimal interpolation problem in Eq. (\ref{eq:intrp}) has been solved earlier using the Conjugate Gradient Method that minimizes the successive \textit{error} at each iteration with respect to the analytical solution, rather than the residuals\cite{kunis2007stability}. Here, the number of iterations must be specified because to measure convergence requires knowledge of an $M\times M$ kernel matrix:
\begin{equation}
    \mathbf{K}_N = A^HW_kA \in \mathbf{C}^{M \times M}
\end{equation}
\noindent which is a circulant matrix\cite{kunis2007stability} that can be evaluated by some $M \times 1$ vector in $\mathcal{O}(MlogM)$ time, but the determination of the matrix itself scales with $M$ in quadratic time. Regardless the the residuals from the ``true'' solution as minimized by the Conjugate Gradient method for the $l^{th}$ iteration can be calculated as:

\begin{align}\label{eq:kernel}
    e_{h_k,l} = ||h_{k,l} - W_kA\mathbf{K}_N^{-1}f(x_j)||_{W_k^{-1}}\\
    = ||h_{k,l} - W_kA\left(A^HW_kA\right)^{-1}f(x_j)||_{W_k^{-1}}.
\end{align}

\begin{figure}
    \centering
    \includegraphics[width=0.95\linewidth]{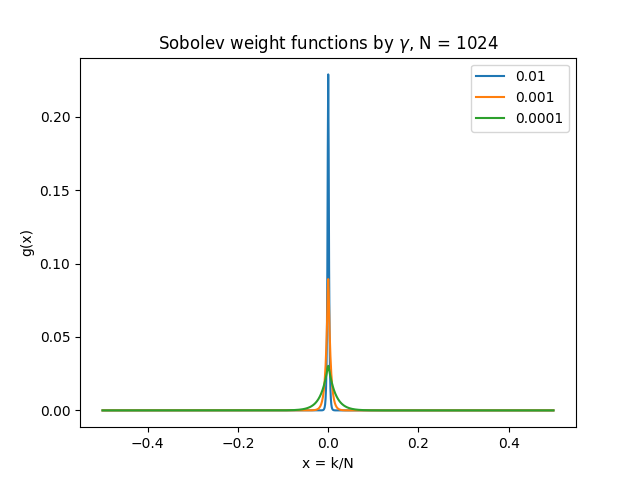}
    \caption{Various values for $\gamma$ and their impact on the Sobolev kernel function. For lower values of $\gamma$, the filter weights higher frequencies more strongly.}
    \label{fig:sobolev}
\end{figure}

Weight functions are used to represent the weight matrix $W_k$ along the diagonal as a low-pass filter to prevent the Fourier coefficients from poorly interpolating discontinuities in the observed data. Admissible weight functions ($V$) are defined by Kunis and Potts  \cite{kunis2007stability} as any function $g \in V$ that satisfies:
\begin{align}
    |g^{\beta-1}|_V = \int_{\mathbf{R}} |dg^{\beta-1}(z)|\\ =sup \sum_{j=0}^{n-1} | g^{\beta-1}(z_{j+1}) - g^{\beta-1}(z_j)| < \infty 
\end{align}
\noindent i.e. that for $|z| < \frac{1}{2}$ the maximum bounded set of the absolute difference between $g^{\beta-1}(z_{k+1})$ and $ g^{\beta-1}(z_k)$ does not diverge when the $\beta^{th}$-fold derivative is taken. The authors further proposed the practical requirement that all admissible weight functions be implicitly normalized by their $L1$ norm to 1. Two possible kernels that satisfy this condition are the modified \textit{Fejer} and \textit{Sobolev} kernels, where the Fejer kernel is a kernel of order $\beta=2$ and of polynomial order $N$:
\begin{equation}
    B_{\beta=2,N}(x) = \frac{2(1 + e^{-2\pi \mathbf{i}x})}{N^2} \left( \frac{\sin(\frac{N}{2}\pi x)}{\sin(\pi x)}\right)^2
\end{equation}
Defining the norm of the \textit{Kernel} function as the convolution of $g * f_j$ via the adjoint transformation of $g$:
\begin{equation}
||g||_{1,N} = \sum_{k= -N/2}^{k=N/2} g\left(\frac{k}{N}\right) e^{2\pi \mathbf{i}kx}
\end{equation}
\noindent can be used to define the Sobolev kernel, first as an intermediate and generalizable weight function:
\begin{equation}\label{eq:weight1}
    g_{\alpha,\beta,\gamma}(z) := c_{\alpha,\beta,\gamma} \frac{\left(\frac{1}{4} - z^2\right)^\beta}{\gamma + |z|^{2\alpha}}
\end{equation}
\noindent which can then be used recursively to solve for the constant: $c_{\alpha,\beta,\gamma}$ by taking the inverse of the sum of the right-most term. The Sobolev kernel is defined as:
\begin{equation}
    S_{\alpha,\beta,\gamma,N} := \frac{1 + e^{-2\pi \mathbf{i} x}}{2 ||g_{\alpha,\beta,\gamma}||_{1,N}} \sum_{k=-N/2}^{N/2}g\left(\frac{k}{N}\right)e^{2\pi \mathbf{i} k x}
\end{equation}
\noindent for some values of $\alpha, \beta, \gamma$. Selection of $\alpha$ and $\beta$ are related to the dimensionality of the transform and conditions for a suitable weight function, but a useful value for $\gamma$ is $10^{-2}$ for 1-dimensional transforms as suggested previously by Kunis and Potts \cite{kunis2007stability}. The resultant weight function favours the low frequency oscillations in the representation of $f(x_j)$ as $h_k$. The relative importance of these frequencies for the Sobolev kernel as a function of $\gamma$ are represented in Figure \ref{fig:sobolev}.

This work demonstrates the quadratic form of Eq. (\ref{eq:intrp}), and offers a convenient method for solving the i-iNFFT without the need to iterate. In addition, it avoids problems with dimensionality by inverting only the $AA^H$ matrix, which sacrifices some time complexity for simplicity. This is an acceptable trade-off for small values of $N$ and is competitive with other algorithms for the direct inversion of non-uniform fast Fourier transforms which require similar $\mathcal{O}(N^2)$ matrix pre-computations\cite{kircheis2019direct}.

\section{Prior Art}\label{sec:lssa}
\subsection{Least squares spectral analyses (LSSA)}
Some proposals exist to solve a similar optimisation problem in Eq.(\ref{eq:min_func}). Van\'{i}\v{c}ek et al. \cite{wells1985least} first proposed a direct inversion to solve for $\hat{h}_k$:
\begin{equation}
    \hat{h}_k = (AA^H)^{-1}Af(x_j)
\end{equation}
which is a simple inversion of the matrix $A$ that executes in $\mathcal{O}(M^2N)$ time, assuming $M > N$. Problems to do with the condition of the matrix can be addressed using any number of rank-deficient solutions to the inverse least squares problem, including the Moore-Penrose pseudo-inverse which for complex matrices is:
\begin{equation}
    A = U_R \Sigma_R V_R^H
\end{equation}
\begin{equation}
    A^+ = V_R \Sigma_R^{-1} U_R^H
\end{equation}
\begin{equation}
    h_k = A^+f(x_j)
\end{equation}
\noindent where $R$ is the rank of the matrix $A$; assumed to be the smallest of either $M$ or $N$ for the best possible reconstruction. 

A more recent implementation calculates the frequency spectrum in $HMlogHM$ time for a \textit{preset} number of sinusoidal components, $H$\cite{palmer2009fast}, via the minimization of the (presumed) Gaussian residuals via the $\chi^2$ statistic:
\begin{equation}
    \chi^2 (\mathbf{P},H) = \sum_{j=1}^H \frac{\left(\mathbf{P}(x_j,H) - f(x_j)\right)^2}{\sigma_j^2}
\end{equation}
\noindent in which the $\chi^2$ distance of the model $\mathbf{P}$ with parameters $H$ is evaluated at each time point $x_j$ relative to the observational data $x_j$ normalized to the apparent variance ($\sigma^2_j$) of the residuals. The $H$  sinusoidal functions are calculated analytically via a linear fit, before being transformed into the frequency domain via a conventional truncated FFT, on the continuously reconstructed function.

The most closely related methods in purpose and time complexity have been proposed by Kircheis et al. \cite{kircheis2019direct} but vary in methodology depending on the conditioning of the matrix. In the case of $M=N$, assuming the transformation matrix of an equidistant time series is $FF^T = MI_N$, a sparse decomposition in the non-equidistant case is sought to approximate an equi-distance transformation matrix, $F$:
\begin{equation}\label{eq:kircheis}
    A = DFB
\end{equation}
\begin{equation}
    MI_N \approx AB^HF^HD^* 
\end{equation}
\noindent Eq. (\ref{eq:kircheis}) contains matrix terms $D$, $F$, and $B$ respectively; they are an $N \times N$ diagonal matrix of entries $N_\sigma \times w_k$, the $N \times M$ truncated Fourier matrix and a $M \times M$ sparse matrix. $N_\sigma$ is the ``oversampling'' factor defined by Potts in \cite{potts2021continuous} for continuous window functions. $w_k$ is a weight function similar to those discussed later in Section \ref{subsec:weights}. The $*$ designation refers to the complex conjugate of the diagonal entries of $D$. The decomposition requires $\mathcal{O}(N^2)$ computational time, but the transform is only $\mathcal{O}(NlogN + M)$ complex \cite{kircheis2019direct}; hence it can be described as a true fast transform when equivalent sparsity is observed in several different calculations. 

\subsection{Related methods for interpolation}

Methods related to spectral power density estimation are not equivalent to least-squares frequency spectra estimation, but reconstruction of irregularly sampled time-series is often a common goal for both. Spectral power density is the amount of \textit{power} transmitted through a signal $f(t)$ integrated across all possible values of time centred around some average time $t_0$:
\begin{equation}
P = \lim_{T\rightarrow \infty} \frac{1}{T} \int_{t_0 - T/2}^{t_0 + T/2} |x(t)|^2 dt.
\end{equation}
When the signal is centred around $x_{t_0}$ this expression simplifies to:
\begin{equation}
P = \lim_{T\rightarrow \infty} \frac{1}{T} \int_{-\infty}^{+\infty} |x_{t_0}(t)|^2 dt,
\end{equation}
\noindent observing that the sum of a squared function transformed as its equivalent Fourier series is equal to the square of the function itself (i.e. Parseval's theorem \cite{tolstov1976fourier}) the definition of power spectral density $P(k)$ for some frequency coefficient $h_k$ is:
\begin{equation}\label{eq:psd}
    P(k) = \lim_{T\rightarrow \infty}\frac{1}{T}|h_k(T)|^2
\end{equation}
\noindent which states that the power of a signal observed over an infinite time series is equal to the square of the Fourier coefficients $h_k$ observed for some time series $T$ centred around $t_0$ \cite{scargle1982studies}. The power density spectrum can be estimated through what is known as a \textit{Periodogram}; the most common method is the \textit{Lomb-Scargle} Periodogram that fits a user-defined number of $N$ periodic functions of the form:
\begin{equation}
    f(x_j) = \sum_{n=1}^N\mathcal{A}_n\sin(\omega_n x_j) + \mathcal{B}_n\cos(\omega_n x_j) + \mu
\end{equation}
\noindent for some $\omega_n$ frequencies and fitted attributes $\mathcal{A}_n$, $\mathcal{B}_n$ along with the overall mean value for the time series $\mu$. The periodical functions are transformed to mutually orthogonal components, and related back to the power spectrum in Eq. (\ref{eq:psd}) through the following relationship that incorporates an offset for each function, $\tau$:
\begin{align}
 P(k) = 
  \frac { \left( \sum_j f(x_j) \cos \omega_n ( x_j - \tau ) \right) ^ 2}
        { \sum_j \cos^2 \omega_n ( x_j - \tau ) } +
\\
 \frac {\left( \sum_j f(x_j) \sin \omega_n ( t_j - \tau ) \right) ^ 2}
        { \sum_j \sin^2 \omega_n ( x_j - \tau ) }
\end{align}
\noindent where $\tau$ is defined as:
\begin{equation}
    \tan 2 \omega_n \tau = \frac{\sum_j \sin 2 \omega_n x_j}{\sum_j \cos 2 \omega_n x_j}.
\end{equation}
From the power spectrum, statistical significance for a limited number of frequencies can be inferred from the $\beta-$distribution and the data reconstructed according to these frequencies\cite{wells1985least}. 

Singular Spectrum Analysis (SSA) is another method that has been used to iteratively impute missing values in time-series data \cite{kondrashov2006spatio}, based on the decomposition of a lag-covariance matrix $C_{a,b}$ calculated on the time series $f(x_j)$:
\begin{equation}\label{eq:covariance}
    C_{a,b} = \frac{1}{M - |a-b|}\sum_{j=1}^{M - |a - b|}f(x_j)f(x_{j + |a - b|}).
\end{equation}
$C_{a,b}$ is an $\Bar{M} \times \Bar{M}$ matrix of the nominal length of the window encompassing the longest possible periodic function. The total lag $|a-b|$ of entries $a$ and $b$ is normalised relative to the total number of populated observations, $M$ minus the relative distance spanned by $|a-b|$; placing less emphasis on regions with more missing data. The resultant matrix is decomposed into its constituent eigenvalues:
\begin{equation}
    C_{a,b}^TC_{a,b} = Q\Lambda Q^{-1}
\end{equation}
\noindent and time series data in $f(t)$ are reconstructed according to the nominal time-indices for the $m^{th}$ eigenvector $m \in [1,\Bar{M}]$:
\begin{equation}
    \hat{f}(t)_m = \sum_{j=1}^{\Bar{M}}f(t + j - 1)Q_m(j).
\end{equation}
The $m^{th}$ oscillatory mode is equivalent to the corresponding Principal Component scores for each time point $t$ projected into the corresponding eigenvectors in $Q$. Similar to a Principal Component Analysis, the relative importance of each oscillatory mode can be inferred from an inspection of the square root of the eigenvalues as a ``scree'' plot of singular values, or any number of other tests for significance can be used to infer their relative importance \cite{vitale2017selecting}\cite{allen1996monte}. For iterative reconstructions of missing data, the reconstructions for some number of principal components are used for the next estimate of the lag-covariance matrix in Eq. (\ref{eq:covariance}) until some convergence criterion is met.

Two main considerations for use of a least-squares solution to the reconstruction of univariate time series data are how sensitive each technique is to the selection of an appropriate number of components (either sinusoids or orthogonal oscillatory reconstructions), and whether or not a convenient closed-form solution exists for its determination. In theory, the use of an appropriate weight function makes the component number selection much less critical, but a weight function cannot be used in instances where a statistical subroutine is used to select a discrete number of frequencies or orthogonal components. Weight functions are common for Fourier analyses, but some problems persist in the non-equispaced cases: either the algorithm for arriving at the solution is scales poorly \cite{wells1985least} with the number of Fourier coefficients, or depends exquisitely on the condition of the transformation matrix \cite{kircheis2019direct}.
\section{Considerations for solving the optimal interpolation problem}\label{sec:solve}

\subsection{Solving for the predicted Fourier coefficients, $\hat{h_k}$}
Let $||h_k||_{W_k^{-1}}$ be the weighted norm of the predicted Fourier coefficients $h_k$, and let $||A^Hh_k - f(x_j)||_F^2$ be the constraint. The optimization problem in Eq. (\ref{eq:intrp}) can be described as:
\begin{equation}\label{eq:cost}
    \argmin_{h_k}\left(||h_k||_{W_k^{-1}} + ||A^Hh_k - f(x_j)||_F^2 \right)
\end{equation}
\noindent which through the norm of each term is equivalent to:
\begin{equation}\label{eq:unpacked}
    \argmin_{h_k}\left(h_k^H W_k^{-1} h_k + \left(A^Hh_k - f(x_j)\right)^H\left(A^Hh_k - f(x_j)\right)\right)
\end{equation}
\noindent and unpacked as a polynomial expression becomes:
\begin{equation}\label{eq:quad}
\begin{split}
           \argmin_{h_k}  (h_k^HW_k^{-1}h_k + h_k^H A A^H h_k\\
   - 2f(x_j)^HA^Hh_k + f(x_j)^Hf(x_j)).
\end{split}
\end{equation}
To find the optimal $h_k$ to minimise Eq. (\ref{eq:unpacked}), let the function $\mathcal{L}$ represent the cost function Eq. (\ref{eq:cost}) for the following derivatives (abbreviating $\mathcal{L} = \mathcal{L}(h_k)$):
\begin{align}
    \frac{\partial \mathcal{L}}{\partial h_k} h_k^H W_k^{-1} h_k = W_k^{-1} h_k + h_k^H W_k^{-1}\\ 
    = 2\hat{W^{-1}}h_k\\
    \frac{\partial \mathcal{L}}{\partial h_k} h_k^H A^H A h_k = 2 AA^Hh_k\\
    \frac{\partial \mathcal{L}}{\partial h_k} h_k^H A f(x_j) = A f(x_j)\\
    \frac{\partial \mathcal{L}}{\partial h_k} f(x_j)^H A^Hh_k = A f(x_j)\\
    \frac{\partial \mathcal{L}}{\partial h_k} f(x_j)^H f(x_j) = 0
\end{align}

Taking $\partial \mathcal{L} / \partial h_k = 0$:
\begin{align} \label{eq:grad}
    2W_k^{-1}h_k + 2 AA^Hh_k - 2 A f(x_j) = 0 \\
    W_k^{-1}h_k + AA^Hh_k = A f(x_j) \\
    h_k = \frac{A f(x_j)}{W^{-1} + AA^H}
\end{align}
where in Eq. (\ref{eq:fun}) the identity $AA^H = MI_N$ holds for equispaced observations as demonstrated in Appendix \ref{sec:app1}. Assuming this identity, even as an approximation in cases were irregular sampling is rare - a solution for $h_k$ can be found by bringing the $W_k^{-1}$ out of the inverse term, for a \textit{weighted, truncated} inverse transform as demonstrated in the following section.

\subsubsection{Optimal interpolation problem for equispaced nodes}

Factoring out the constant $M$:
\begin{equation}\label{eq:fun}
    h_k = \frac{Af(x_j)}{M\left(\frac{1}{M}\hat{W}^{-1}_k + I\right)}
\end{equation}
\noindent it is possible to use one of the Searle identities\cite{petersen2008matrix}: $(I + X^{-1})^{-1} = X(X + I)^{-1}$ so that Eq. (\ref{eq:fun}) becomes the simplified expression for the gradient as a function of the weights:
\begin{equation}
    \nabla h_k = \frac{Af(x_j)W_k}{MW_k + I_N}
\end{equation}

\noindent since $\frac{1}{M}W_k^{-1} = (MW_k)^{-1}$ for a diagonal matrix.

\subsubsection{Optimal interpolation problem for non-equispaced nodes}

When $f_j \in [-\frac{1}{2}, \frac{1}{2})$ is not equispaced, the identity: $AA^H = MI_N$ does not hold (Appendix \ref{sec:app1}). The most straightforward way of solving for $h_k$ in Eq. (\ref{eq:grad}) is via a numerical approximation of $AA^H$; however warranting further consideration is the fact that it exists as part of a total matrix inverse as: $(W_k^{-1} + AA^H)^{-1}$. A direct inverse of this expression is not numerically stable, since $\lim_{n \rightarrow |N/2|}W_{k=n} = 0$, $\lim_{n \rightarrow |N/2|}W_{k=n}^{-1} = \infty$ for the weighted Fourier coefficients representing the higher frequencies. We are motivated to find a more stable numerical representation of  Eq. (\ref{eq:grad}) via the Kailath identity \cite{bishop1995neural,petersen2008matrix}:
\begin{equation}\label{eq:kail}
    \left(X + YZ\right)^{-1} = X^{-1} - X^{-1}Y\left(I + ZX^{-1}Y\right)^{-1}ZX^{-1}
\end{equation}
\noindent in which we can substitute $Y,Z$ in Eq. (\ref{eq:kail}) via the upper and lower triangular matrices of $AA^H$ via the $LU$ decomposition:
\begin{align}\label{eq:subkail}
    AA^H = LU\\
    \left(W_k^{-1} + LU\right)^{-1} = W_k - W_kL\left(I + UW_kL\right)^{-1}UW_k.
\end{align}
Substituting Eq. (\ref{eq:subkail}) into the denominator of  (\ref{eq:grad}), the function that minimises the cost function for the optimal interpolation problem as a function of $h_k$ in  Eq. (\ref{eq:cost}):
\begin{align}\label{eq:final}
    \hat{h}_k = A f(x_j) \left(W_k - W_kL\left(I + UW_kL\right)^{-1}UW_k\right)
\end{align}
\noindent where it is clear that in Eq. (\ref{eq:final}) we avoid the problem of directly inverting the kernel function in the inverse term of the cost function along the diagonal of $W_k$, as a substitute for Eq. (\ref{eq:grad}).

\subsection{Second derivative test}
From Eq. (\ref{eq:grad}), the second derivatives with respect to $h_k$ are as follows:
\begin{align}
    \frac{\partial ^2 \mathcal{L}}{\partial h_k^2} h_k^H W_k^{-1} h_k = 2W_k^{-1}\\
    \frac{\partial^2 \mathcal{L}}{\partial h_k^2} h_k^H A A^H h_k = 2AA^H
\end{align}
\noindent and again, since the inverse of any kernel function $\lim_{n \rightarrow |N/2|}W_{k=n} = 0$, $\lim_{n \rightarrow |N/2|}W_{k=n}^{-1} = \infty$, we are motivated to reduce the form:
\begin{align}
    \nabla^2 h_k = AA^HW_k + I_N.
\end{align}
$AA^H$ (as a symmetric matrix for which $rank(A) = N$) and $W_k$ (as a positive, symmetrical, diagonal matrix) are positive semi-definite (PSD) matrices, so their product must also be PSD \cite{meenakshi1999product}. The identity matrix $I_N$ does not change this property. Consequently, the second derivative test indicates that the step direction \textit{minimises} Eq. (\ref{eq:unpacked}) and that the cost function is convex.

\section{Analysis of Remote Sensor Data}\label{sec:res}

\subsubsection{Remote Sensor Data}
Soil temperature data were collected in the Sierra Nevada National Park in Andaluc\'ia, Spain, as part of the Smart EcoMountains program in partnership with the Lifewatch ERIC consortium. Solar-powered, remote, and multi-parametric micro-stations were placed within different environments (shrub, rock, and oak litter) and altitudes (from 1300m to 2300m), for the dates inclusive of 29-07-2020 until 04-05-2023 when the data were queried for analysis. Measurements were continuously collected at a nominal 10 minute interval with sensors (ibuttons, thermochron-type; Maxim integrated) connected to a network that transmits the collected information when connectivity is available. All data used in this study were visually inspected for extreme values outside expected ranges that could occur in the time series due to sensor failures. These erroneous values were reported as NaNs.

\subsubsection{Computational information}
Calculations were run with a B550 AORUS ELITE AX V2 motherboard with 128 GiB DDR4 RAM, an AMD Ryzen 9 5950x 32 core processor with an NVIDIA GeForce RTX 3080 graphics processor, and a 2 TB NVMe Samsung SSD 980. Software was installed on Ubuntu 22.04.3 LTS 64-bit operating system ``Jammy Jellyfish''.

Graphical inspections and data cleaning were conducted with the tidyverse 2.0.0 package \cite{tidyverse} implemented in the R 4.3.1 software environment \cite{R_Software}.

Numerical routines were performed within the Python 3.10 release. Data were parsed from the cleaned .csv values using Pandas 2.0.3. Various numerical operations were performed using the vectorized functions available in NumPy 1.25 \cite{harris2020array}, with additional statistical functionality from SciPy 1.11 \cite{2020SciPy-NMeth}. Graphics were created using Matplotlib 3.7.1 \cite{Hunter:2007}. 

\begin{figure}[hbt!]
    \centering
    \includegraphics[width=0.95\linewidth]{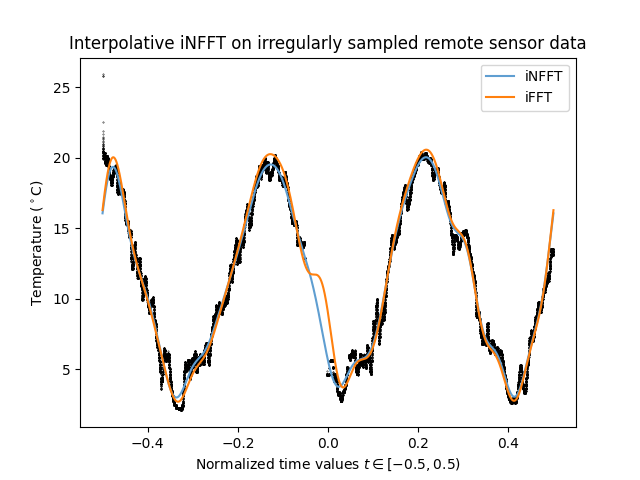}
    \caption{Results of the i-iNFFT algorithm using the Sobolev kernel with $\gamma = 1e-2$, and $N=1024$. The mean absolute fractional error: $\frac{1}{M}\sum_{j=0}^M|(y_{pred} - y_{obs}) / y_{obs}|$ is $6.76 \times 10^{-2}$ for the i-iNFFT reconstruction, versus $7.48 \times 10^{-2}$ for the equivalent weighted truncated inverse FFT (t-iFFT) reconstruction. Note the stronger tendency of the t-iFFT to regress towards the mean for larger gaps.}
    \label{fig:interp}
\end{figure}

In order to validate the solution using previously unconsidered data, measurements were randomly removed from a largely continuous time series. The i-iNFFT was trained on the remaining data, before the data that was previously removed was used as a test set to evaluate the agreement of the interpolated function with the previously unconsidered data. In the second test, a contiguous block of data was removed and a similar procedure was followed. A visual example of the interpolation performance over a missing contiguous block can be seen in Figure \ref{fig:interp}, where a large number values failed a QC check over a period of several months. These tests were performed using data with agnostic, or na\"ieve time labels, in the sense that the earliest recorded observations were indicated as $-0.5$, and the last recorded observations were indicated at just under $0.5$. This does not reflect the true periodicity of the data, and while it could be inferred from the time labels themselves, these labels were not applied on these data for the sake of demonstrating both the windowing effect towards the beginning and ends of the time series, and the potential to infer the periodicity of the time series through stochastic gradient descent. 

On any one external set, the performance of the i-iNFFT algorithm depends on the $N$ number of Fourier coefficients, and the value of $\gamma$. The results suggest that for an appropriate value of $\gamma$ as mentioned in previous literature \cite{kunis2007stability}, the results are only weakly dependent on the $N$ number of Fourier coefficients.

\subsection{Performance metrics relative to truncated Fourier transforms}

All training data were mean-centred prior to analysis and these means were added to the reconstructed data prior to recording the residuals. The results summarize cross validation routines for randomly sampled data in Tables \ref{tab:corr_rand} and \ref{tab:predrand}, and randomly sampled contiguous blocks in Tables \ref{tab:corr_blck} and \ref{tab:predblck} at 3 levels of missing data (10, 20, and 30\%) sampled 20 times to get the average error, and the uncertainty, as one standard deviation. Both relative error for the predicted set, and the correlation coefficient for the reconstructed data on the full data set versus the partially sampled dataset were also calculated to intuit the error due to the constraint of the reconstruction (via the weighted kernel approximation) versus the irregularity in the sampling of the observed data. The results of each experiment suggest that the i-iNFFT algorithm outperforms the equivalent inverse Fast Fourier Transform (iFFT) algorithm, as demonstrated by a one-sided permutational t-test. This demonstrates that in addition to the mean of the performance statistic being higher in all cases, this advantage is maintained consistently throughout all possible combinations of the cross-validation results while not assuming an underlying normal distribution of the residuals about the mean. 

\begin{table}[hbt!]
\begin{center}
\begin{tabular}{lll}
                        & $r^{i-iNFFT}_{|rand}$                       & $r^{t-iFFT}_{|rand}$                        \\ \cline{2-3} 
\multicolumn{1}{l|}{10} & \multicolumn{1}{l|}{1.00E+00 $\pm$ 3E-07 **} & \multicolumn{1}{l|}{9.99E-01 $\pm$ 4E-05} \\ \cline{2-3} 
\multicolumn{1}{l|}{20} & \multicolumn{1}{l|}{1.00E+00 $\pm$ 9E-07 **} & \multicolumn{1}{l|}{9.99E-01 $\pm$ 6E-05} \\ \cline{2-3} 
\multicolumn{1}{l|}{30} & \multicolumn{1}{l|}{1.00E+00 $\pm$ 2E-06 **} & \multicolumn{1}{l|}{9.99E-01 $\pm$ 8E-05} \\ \cline{2-3} 
\end{tabular}
\caption{Comparison of the i-iNFFT, and t-iFFT reconstruction errors based on correlation between the trigonometric polynomial trained on all available data, versus the randomly sampled data at 10, 20, and 30\% external sampling.** denotes a permutational, one-sided t-test indicating $p<0.01$ ($N=20, \mathcal{\pi}$ (permutations) $=10,000$)}
\label{tab:corr_rand}
\end{center}
\end{table}
\begin{table}[hbt!]
\begin{center}
\begin{tabular}{lll}
                        & $r^{i-iNFFT}_{|blck}$                       & $r^{t-iFFT}_{|blck}$                        \\ \cline{2-3} 
\multicolumn{1}{l|}{10} & \multicolumn{1}{l|}{9.91E-01 $\pm$ 6E-03 **} & \multicolumn{1}{l|}{9.92E-01 $\pm$ 4E-02} \\ \cline{2-3} 
\multicolumn{1}{l|}{20} & \multicolumn{1}{l|}{9.61E-01 $\pm$ 3E-02 **} & \multicolumn{1}{l|}{8.72E-01 $\pm$ 5E-02} \\ \cline{2-3} 
\multicolumn{1}{l|}{30} & \multicolumn{1}{l|}{8.83E0-1 $\pm$ 7E-02 **} & \multicolumn{1}{l|}{7.97E-01 $\pm$ 4E-02} \\ \cline{2-3} 
\end{tabular}
\caption{Comparison of the i-iNFFT, and t-iFFT reconstruction errors based on correlation between the trigonometric polynomial trained on all available data, versus the contiguous block sampled data at 10, 20, and 30\% external sampling. ** denotes a permutational, one-sided t-test indicating $p<0.01$ ($N=20, \pi=10,000$)}
\label{tab:corr_blck}
\end{center}
\end{table}
\begin{table}[hbt!]
\begin{center}
\begin{tabular}{lll}
                        & $(1-Err)^{i-iNFFT}_{|rand}$                 & $(1-Err)^{t-iFFT}_{|rand}$                  \\ \cline{2-3} 
\multicolumn{1}{l|}{10} & \multicolumn{1}{l|}{9.86E-01 $\pm$ 3E-04} & \multicolumn{1}{l|}{9.86E-01 $\pm$ 3E-02} \\ \cline{2-3} 
\multicolumn{1}{l|}{20} & \multicolumn{1}{l|}{9.73E-01 $\pm$ 5E-04} & \multicolumn{1}{l|}{9.73E-01 $\pm$ 6E-02} \\ \cline{2-3} 
\multicolumn{1}{l|}{30} & \multicolumn{1}{l|}{9.59E-01 $\pm$ 6E-04} & \multicolumn{1}{l|}{9.59E-01 $\pm$ 7E-02} \\ \cline{2-3} 
\end{tabular}
\end{center}
\caption{Comparison of the i-iNFFT, and t-iFFT reconstruction errors based on external prediction for randomly selected data at 10, 20, and 30\% external sampling. The i-iNFFT did not demonstrate a statistically significant improvement in this case according to a permutational one-sided t-test ($N=20, \pi =10,000$)}
\label{tab:predrand}
\end{table}
\begin{table}[hbt!]
\begin{center}
\begin{tabular}{lll}
                        & $(1-Err)^{i-iNFFT}_{|blck}$                 & $(1-Err)^{t-iFFT}_{|blck}$                  \\ \cline{2-3} 
\multicolumn{1}{l|}{10} & \multicolumn{1}{l|}{9.91E-01 $\pm$ 6E-03 *} & \multicolumn{1}{l|}{9.28E-01 $\pm$ 5E-02} \\ \cline{2-3} 
\multicolumn{1}{l|}{20} & \multicolumn{1}{l|}{9.61E-01 $\pm$ 3E-02 *} & \multicolumn{1}{l|}{8.67E-01 $\pm$ 7E-02} \\ \cline{2-3} 
\multicolumn{1}{l|}{30} & \multicolumn{1}{l|}{8.83E-01 $\pm$ 7E-02} & \multicolumn{1}{l|}{8.12E-01 $\pm$ 4E-02} \\ \cline{2-3} 
\end{tabular}
\end{center}
\caption{Comparison of the i-iNFFT, and t-iFFT reconstruction errors based on external prediction for contiguous blocks of selected data at 10, 20, and 30\% external sampling. * denoted a permutational, one-sided t-test indicating $p<0.05$ ($N=20, \pi=10,000$)}
\label{tab:predblck}
\end{table}

The relative error for the randomly sampled data is on the same order of magnitude as reported in the original study for the iterative solution to the optimal interpolation problem. However, it is clear that a more critical metric for determining the accuracy of the interpolation is by removing contiguous blocks of data, and comparing the reconstruction against the true values for those data. Despite this, and considering that the interpolation was ``smoothed'' using a rather aggressive low-pass filter, the results are reasonable and corroborate that the calculated spectrum of the data is relatively consistent with what data is observed, which is what close to what we would expect from a true inverse transform.

\subsection{Effect of $\gamma$ and $N$ on reconstruction}

Using a single time series collected as a part of the original study, but different than the one analyzed previously, the following experiments were performed to assess the impact of the selection of the $\gamma$ coefficient and the number of Fourier coefficients on the relative reconstruction accuracy on both the known values and a randomly selected, contiguous external prediction set ($20\%$ of the total number of observations in the dataset).

\begin{figure}[hbt]
    \centering
    \includegraphics[width=0.95\linewidth]{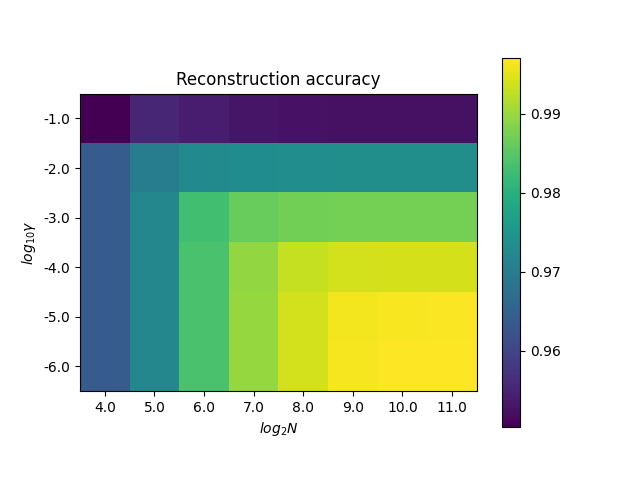}
    \caption{Reconstruction of irregularly sampled time-domain data by $\gamma \in [1e-1,1e-6]$ and $N \in [2^{4},2^{11}]$}
    \label{fig:int_coeff}
\end{figure}

\begin{figure}[hbt]
    \centering
    \includegraphics[width=0.95\linewidth]{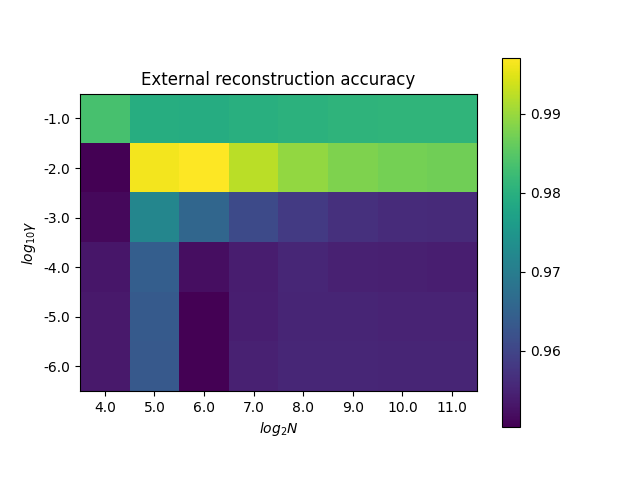}
    \caption{Reconstruction of irregularly sampled time-domain data by $\gamma \in [1e-1,1e-6]$ and $N \in [2^{4},2^{11}]$ with 20\% missing data as one contiguous block}
    \label{fig:ext_coeff}
\end{figure}

The reconstruction accuracy relative to the known data increases as $\gamma$ is lowered over several orders of magnitude, and as the number of Fourier coefficients (here as powers of 2) increases. As far as the Fourier coefficients are concerned, this is an unsurprising result, but it illustrates the effect of $\gamma$ on the reproduction of the original data as displayed in Figure \ref{fig:int_coeff}. 

On an external set, the performance displays an obvious dependence on the values for $N$ and for $\gamma$ but in a way that is not immediately obvious (Figure \ref{fig:ext_coeff}). Certain combinations of these two factors can impact the predictive ability of the Fourier model, but in general increasing the number of coefficients decreases its interpolative ability; as does decreasing the value of the $\gamma$ coefficient. Further discussions on windowing functions for the NFFT have been elaborated on in further work by Potts et al. \cite{potts2021continuous}, are are deserving of further consideration as they relate to this work. 

\subsection{Theoretical time complexity}

In Eq. (\ref{eq:final}), the equation can be broken into the weighted inverse of the $AA^H$ term, which can be solved in theoretical $\mathcal{O}(N^2)$ time and the $Af(x_j)$ term which can be solved in (approximately) $\mathcal{O}(MlogM)$ time as part of an NFFT step. The theoretical time complexity of the weighted inverse term is derived from the fact that Toeplitz matrices can be decomposed via an $LU$ decomposition in $\mathcal{O}(N^2)$ time, and the multiplication of either the $L$ or $U$ term as upper and lower sparse triangular matrices respectively, scale similarly. The current implementation sacrifices portability for time complexity, and performs the $LU$ decomposition in $\mathcal{O}(MN^2)$ complexity via the inner product of the non-uniform transformation matrix, $A$. Although with knowledge of the number of missing values, this function could be automated using a lower-level routine.
The previous algorithm for the inverse interpolative NFFT solves the problem iteratively in $\mathcal{O}(MlogM)$ via the complexity of the adjoint NFFT in which the number of observations in the time domain $M$, is the dominant term.

\subsection{Problem of na\"{i}eve time labels}\label{subsec:time_labels}

In the case of data with labelled time and known periodicity, it is possible to infer the beginning and end of the most significant oscillations within the data as $x_j$. However this knowledge is not always known \textit{a-priori}, and even in cases where it is known such assumptions will not always be fulfilled. In the case of climate sensor data, it cannot necessarily be assumed that an annual recurrence of a temperature trend begins and ends on the same date each year. Incorrectly assuming that the observed data in $x_j$ exhibits a periodicity that starts and ends at the same time each year may result in an effect towards the beginning and end of the observed time-series as shown in Figure \ref{fig:interp}.

\begin{figure}[hbt]
    \centering
    \includegraphics[width=0.95\linewidth]{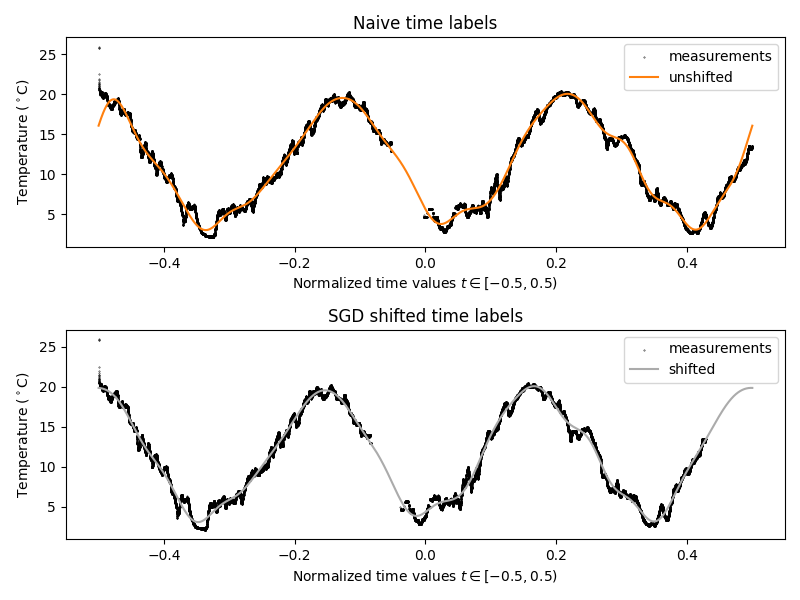}
    \caption{Calculating optimal time labels using a variation of the SGD algorithm. Results using the na\"ieve time labels are shown above, and the calculated time labels are shown below. The error in the reconstruction of the data with na\"ieve time labels was $6.76 \times10^{-2}$, versus $1.41 \times 10^{-2}$ for the calculated time labels via the absolute mean fractional error. }
    \label{fig:sgd}
\end{figure}

Fixing the beginning of a cycle to be the time point as: $-0.5$, and an arbitrary point: $-0.5 < m_e < 0.5 - 1/M$ as encompassing a total fraction of an irregularly sampled period function $f(x_j)$, we can measure the effect of a shift in the observed phase of the data as it affects its reconstruction. The new, nominal number of regularly-sampled observations can be calculated from $m_e$ and $M$:

\begin{equation}
    M_{ext} = \frac{M-1}{m_e + 0.5}
\end{equation}

\noindent and in cases where $m_e$ is a better estimate, we expect the tailing effect to decrease and likewise the total residuals of the interpolated time-series function relative to the known data. While somewhat computationally intensive, this can be optimised using a variant of stochastic gradient descent (SGD).
\begin{algorithm}[hbt]
    \SetAlgoLined
    \KwIn{Initial stop index $m_e$, learning rate $\alpha$, number of batches: $\mathcal{B}$, number of epochs: $\mathcal{E}$}
    \KwOut{Optimal stop index: $m_e$ that minimises the cost function $E(m_e)$}
    \For{$i \leftarrow 1$ \KwTo $\mathcal{E}$}{
        \For{$j \leftarrow 1$ \KwTo $\mathcal{B}$}{
            Compute gradient: $\nabla E(m_e) \leftarrow \argmin_{m_{e,i}} E(m_{e,i})- E(m_{e,j})$\;
            \If{$E(m_{e,i}) < E(m_{e,j})$}{
            Update parameters: $m_e \leftarrow m_e - \alpha \cdot \nabla E(m_e)$\;
            $E(m_{e,j}) \leftarrow E(m_{e,i})$
            }
        }
    }
    \caption{SGD Algorithm for optimizing the stop index, $m_e$, of the observed time series data.}
    \label{alg:sgd}
\end{algorithm}
Since the index $m_e$ is selected randomly at each iteration, several batches are necessary to estimate an informative gradient. In this example, a total of 7 epochs were used with 12 batches for each epoch that were computed in parallel using the \textit{multiprocessing} module in Python. The results for this analysis, using Algorithm \ref{alg:sgd} are shown in Figure \ref{fig:sgd}.
At an optimum value for $m_e$, it can be assumed that the periodicity of the data is some multiple of observed time points. In Figure \ref{fig:sgd} there are three observed annual oscillations, which is consistent with the timestamps on the original data.

\section{Conclusion}
As a whole, the results suggest that the interpolative algorithm performs significantly better relative to the equivalent truncated inverse Fast Fourier Transform, or t-iFFT using an equivalent number of coefficients and weight function in virtually all cases. The differences are most differentiated relative to correlation with the ``true'' underlying trigonometric polynomial. Despite poorer evidence of significance in one of the four tested cross-validation routines, on average the i-iNFFT performs better than the t-iFFT for up to 30\% missing data. There is a significant ``windowing'' effect at the edges of the time series when the relative position of the cycles was not inferred, and this affects the results of the analysis. Nonetheless, this error is consistent in order of magnitude with the original publication of the interpolative transform and can be improved using either \textit{a-priori} knowledge of the data, or using the SGD algorithm demonstrated in this work. The proposed formulation will find use in a variety of signal processing applications, where the use of an FFT is impractical given its assumption of a consistent time-domain density. The software package used in this work is available online: \url{https://github.com/CoDaSLab/intrp_infft_1d}

From an ecological perspective, the findings represent an important step in the handling of climate data acquired through remote sensors, that enable the precise characterization of the micro-climatic conditions within high mountain ecosystems and the organisms that inhabit them. The improved understanding of micro-climatic nuances will bridge the gap between commonly available macro-climatic data and local topographical complexities, enabling the design of more effective conservation strategies and adaptation measures to address the impacts of climate change in these remote and ecologically sensitive regions.

The problem of irregular sampling in high resolution time-series analysis is becoming more common where large volumes of data are required to interrogate hypotheses of complex systems. Frequency-domain analysis with conventional FFTs can introduce considerable errors where an average sampling rate cannot be assumed, but here a convenient mathematical expression to solve for the spectrum of irregularly sampled time-domain data has been demonstrated in a way that satisfies the conditions for an inverse non-uniform Fourier transform. This methodology has been derived from the formalization of the problem proposed earlier, and via the quadratic form of the proposed cost function we can ascertain that the solution is unique and convex. Using a randomly sampled time series from data with few missing values, the observed relative error is on the same order of magnitude as the most current solution to this problem \cite{kunis2007stability}. Using cross-validation with contiguous blocks, we report a relative error that is much higher, but this is a more critical lens through which to view equivalent problems. This is demonstrated very well by the performance of the proposed i-iNFFT algorithm relative to its equivalent, weighted t-iFFT. 

\appendices

\section{Proof of hermitian self-adjoint product identity for equidistant time-domain measurements}\label{sec:app1}
Consider the case of an NFFT matrix and its hermitian, where:
\begin{align}\label{eq:nfftmat}
    A = \left(e^{-2\pi\mathbf{i} k_f \frac{j}{M}}\right)_{k_f \in N,j \in M}\\
    A^H = \left(e^{2\pi\mathbf{i} k_a \frac{j}{M}}\right)_{k_a \in N,j \in M } 
\end{align}
\noindent and $j = -\frac{M}{2}, ..., \frac{M}{2} - 1$, $k_f$ (representing the frequency-wise columns of the forward transform) $= k_a$ (representing the same columns for the adjoint transform) $= -\frac{N}{2}, ..., \frac{N}{2} - 1$. Multiplying the $N \times M$ matrix $A$ by its adjoint yields the $N \times N$ matrix:
\begin{equation}\label{eq:inner}
    AA^H = \sum_{j=-M/2}^{M/2-1} e^{-2\pi\mathbf{i} j k_f} e^{2\pi\mathbf{i} k_a j}
\end{equation}
Eq. (\ref{eq:inner}) simplifies to:
\begin{equation}\label{eq:zero}
    AA^H = \sum_{j=-M/2}^{M/2-1} e^{2\pi\mathbf{i} j (k_a-k_f)}.
\end{equation}
It is clear from Eq. (\ref{eq:zero}) that $e^{2\pi\mathbf{i} j(0)} = 1$ where $k_f=k_a$. This indicates that the product of $AA^H$ along the diagonal is $M$ for the summation of $j \in M$ in Eq. (\ref{eq:zero}). To prove that the off-diagonal matrices are zero, consider $\theta = 2\pi(k_a-k_f)$ as a complex number, $z = \cos{\theta} + \mathbf{i}\sin{\theta}$ following Euler's identity:
\begin{equation} \label{eq:Euler}
    \sum_{j=-M/2}^{M/2-1}e^{\mathbf{i}\theta j} = z
\end{equation}
Multiplying both sides of Eq. (\ref{eq:Euler}) by $e^{\mathbf{i}\theta}$ yields the following expression:
\begin{equation} \label{eq:arbitrary_eith}
    ze^{\mathbf{i}\theta} = \sum_{j=-M/2}^{M/2-1}e^{\mathbf{i}\theta(j+1)}\\
\end{equation}
\noindent since $\sum_{j=0}^{M-1}e^{\mathbf{i}\theta(j+1)}$ is a $j^{th}$ integer multiple of $\theta$, $\sum_{j=0}^{M-1}e^{\mathbf{i}\theta(j+1)} = \sum_{j=0}^{M-1}e^{\mathbf{i}\theta j}$. It follows that:
\begin{equation} \label{eq:end_of_proof}
    ze^{\mathbf{i}\theta j} = z.
\end{equation}
Evaluating Eq. (\ref{eq:end_of_proof}), either $\theta = 0$ in the case of Eq. (\ref{eq:zero}), or for instances in which $k_f \neq k_a$, $z = 0$. This proves that the off-diagonal entries of $AA^H$ are zero, proving the identity for instances where the $j^{th}$ time domain observations are equidistant:
\begin{equation}
    AA^H = M I_N
\end{equation}
\section{$AA^H \neq MI_N$ in the non-equidistant case}\label{sec:app2}

Consider again the equispaced example, where $j \in [-\frac{M}{2},\frac{M}{2})$. We seek to prove that for some complex number, $z$ evaluated as the summation of $\sum_{j = -M/2}^{M/2 -1} e^{\mathbf{i}\theta j} = z$, $z = 1$ when $\theta = 0$ and $z = 0$ when $\theta \neq 0$. We examine the same identity, but in this case using an evaluation of the finite geometric series, where:
\begin{equation}
    \sum_{i=n}^{n-1} r^k = \frac{1-r^n}{1-r}. 
\end{equation}
\noindent since it is obvious that when $\theta = 0$, the only possible evaluation for $z$ is $1$, consider only when $\theta \neq 0$, which we can write as:
\begin{equation}
    z = e^{-\mathbf{i}\theta \frac{M}{2}} + \sum_{j = -M/2 + 1}^{M/2 -1} e^{\mathbf{i}\theta j}
\end{equation}
\noindent to evaluate for $z$ for some constant integer $\theta$, consider the finite geometric series:
\begin{equation}\label{eq:quant1}
    \sum_{j = -M/2}^{M/2 -1} e^{\mathbf{i}\theta j} = \frac{1 - e^{\mathbf{i}\theta (M - 1)}}{1-e} = 0.
\end{equation}
Because $\sum_{-\frac{M}{2} + 1}^{\frac{M}{2} - 1}$ does not evaluate to a convenient integer, evaluating it terms of known quantities:
\begin{equation}\label{eq:solvefor}
    \sum_{j = -M/2}^{M/2 -1} e^{\mathbf{i}\theta j} = e^{-\mathbf{i}\theta \frac{M}{2}} + \sum_{j = -M/2 + 1}^{M/2 -1} e^{\mathbf{i}\theta j}
\end{equation}
\noindent since $e^{-\mathbf{i}\theta \frac{M}{2}} = - 1$ and the LHS of Eq. (\ref{eq:solvefor}) was evaluated in Eq. (\ref{eq:quant1}), solving for $\sum_{j = -M/2 + 1}^{M/2 -1} e^{\mathbf{i}\theta j}$ yields $1$. Therefore, for some complex number $z$, when $\theta \neq 0$, $z = 0$.

Using this alternative formula, it is easy to prove that in the non-equidistant case where $x_j \in [-\frac{M}{2}, \frac{M}{2} - 1)$ the geometric series is not guaranteed, or even likely to sum to $M-1$. And so from Eq. (\ref{eq:quant1}):
\begin{equation}
\forall(\theta \neq 0), \exists \left( z : z = \sum_{j \in x_j} e^{\mathbf{i}\theta j} \neq 0 \right)
\end{equation}

The proof for the equidistant case was adapted from a proof by David J. Fleet and Allan D. Jepson from the University of Toronto \cite{fourier}. The authors would like to acknowledge support from TED Project TED2021-130888B-I00: Biogenic Refuges as modulators of Climate Change in mountain ecosystems (Mountain BIOREFUGES) for this project.

\ifCLASSOPTIONcaptionsoff
  \newpage
\fi

\bibliographystyle{IEEEtran}
\bibliography{main}
\begin{IEEEbiography}[{\includegraphics[width=2.5cm]{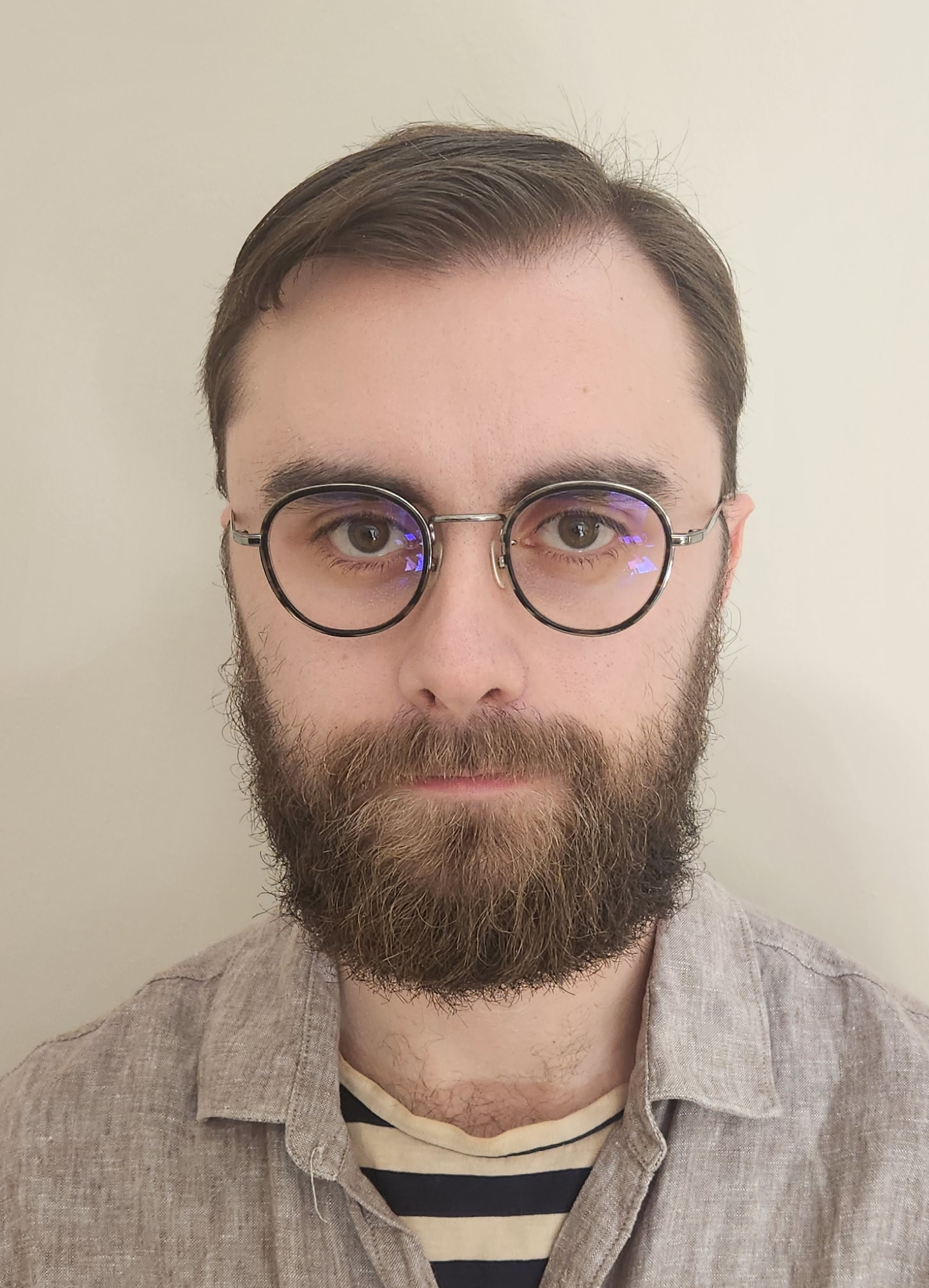}}]{Michael Sorochan Armstrong}
Michael Sorochan Armstrong is a Postdoctoral Research Fellow with the Computational Data Science (CoDaS) Lab in the Department of Signal Theory, Telematics and Communications in the University of Granada, Spain. He completed his PhD in Analytical Chemistry from the University of Alberta, Canada in 2021. His research interests include numerical optimization, and interpretable machine learning. 
\end{IEEEbiography}
\begin{IEEEbiography}[{\includegraphics[width=2.5cm]{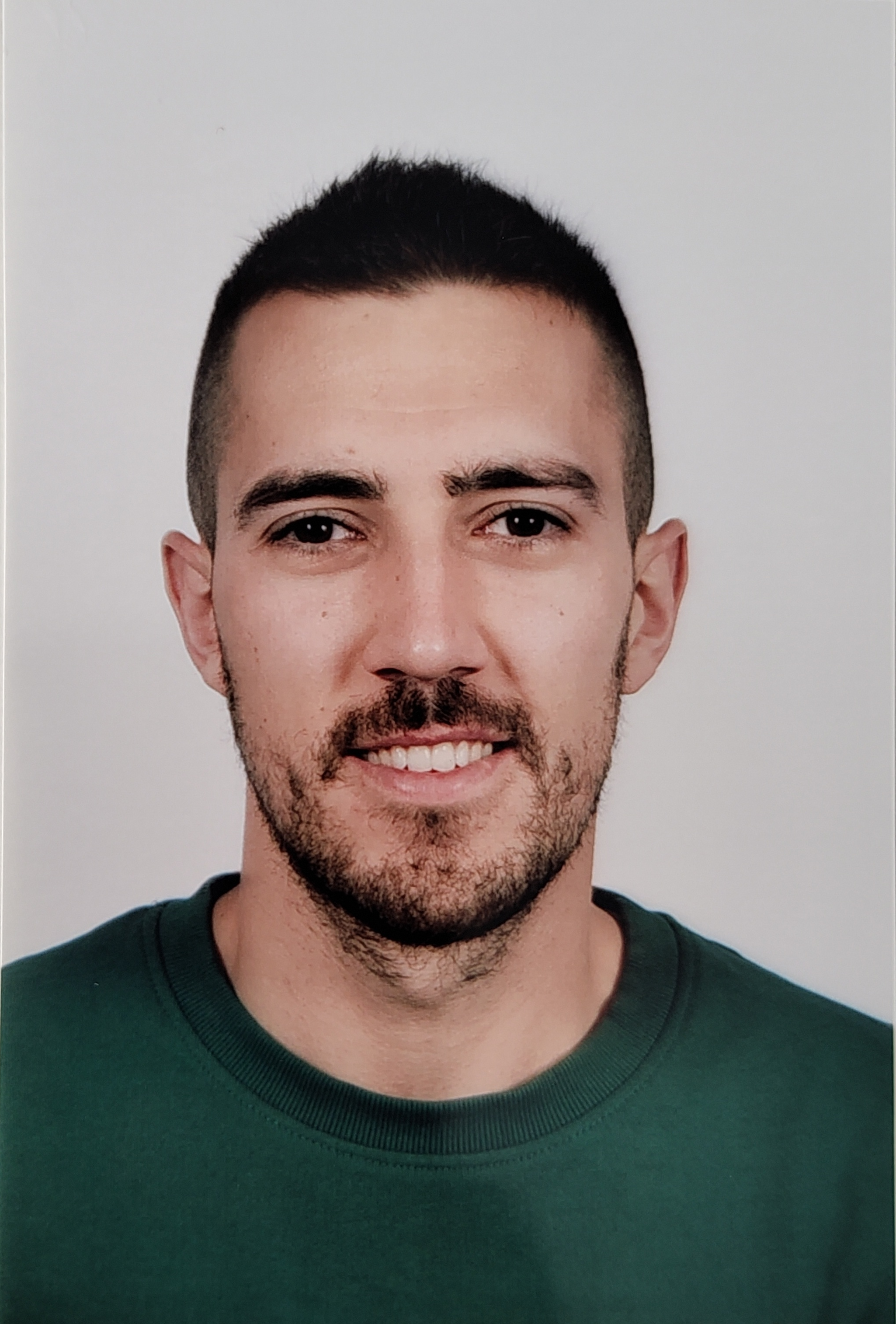}}]{José Carlos Pérez-Girón}
José Carlos Pérez-Girón is a Forestry Engineer (2016) and a Ph.D. in Biogeosciences from the University of Oviedo (2022), granted by the Severo Ochoa program. He serves as a Postdoctoral researcher in the Department of Botany and in the Inter-university Institute for Research on the Earth System in Andalucia through the University of Granada, Spain. His research interests include climate change, ecological modelling and ecosystem services.
\end{IEEEbiography}
\begin{IEEEbiography}[{\includegraphics[width=2.5cm]{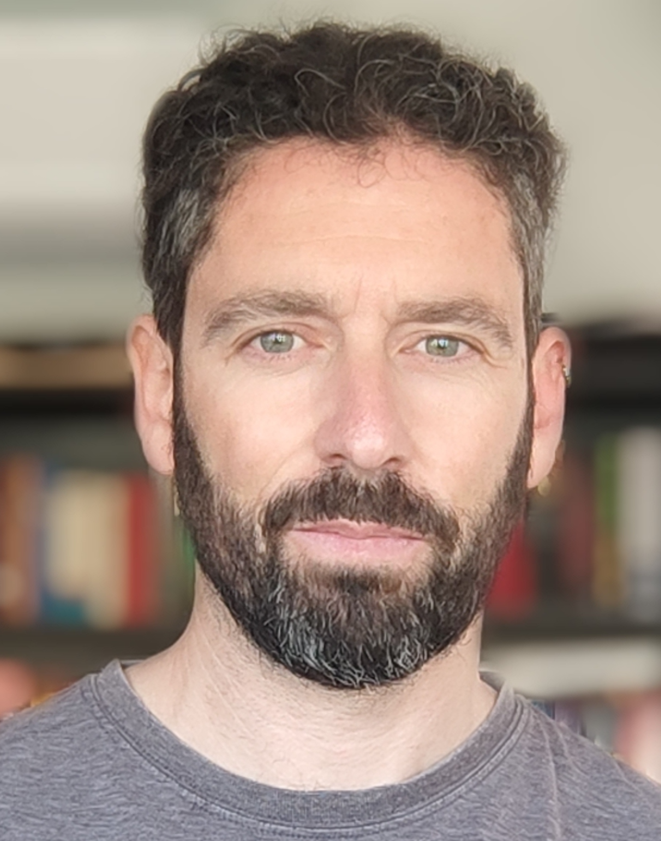}}]{José Camacho}
José Camacho is Full Professor in the Department of Signal Theory, Telematics and Communication and head of the Computational Data Science Laboratory (CoDaS Lab), at the University of Granada, Spain. He holds a degree in Computer Science from the University of Granada (2003) and a Ph.D. from the Technical University of Valencia (2007), both in Spain. He worked as a post-doctoral fellow at the University of Girona, granted by the Juan de la Cierva program, and was a Fulbright fellow in 2018 at  Dartmouth College, USA. His research interests include networkmetrics and intelligent communication systems, computational biology, knowledge discovery in Big Data and the development of new machine learning and statistical tools.
\end{IEEEbiography}
\begin{IEEEbiography}[{\includegraphics[width=2.5cm]{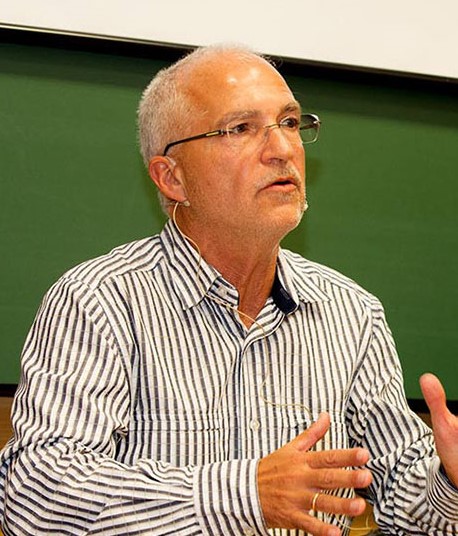}}]{Regino Zamora} Regino Zamora is a Professor of Ecology at the University of Granada. His main field of research is the study of ecological interactions and their consequences on communities and ecosystem processes under global change scenarios. Its academic and research activity is complemented by ongoing collaboration with the public administrations responsible for environmental management, promoting the transfer and application of scientific knowledge to the conservation, management and restoration of ecosystems.  
\end{IEEEbiography}







\end{document}